\documentclass{article}

\begin{document}

\title{Two-forms on four-manifolds and elliptic equations}
\author{S. K. Donaldson}
\maketitle



\newtheorem{thm}{Theorem}
\newtheorem{question}{Question}
\newtheorem{prop}{Proposition}
\newtheorem{defn}{Definition}
\newtheorem{conj}{Conjecture}
\newtheorem{condn}{Condition}
\newtheorem{cor}{Corollary}
\newtheorem{lem}{Lemma}

\section{Background}

Let $V$ be a four-dimensional real vector space with a fixed orientation.
Then the wedge product can be viewed, up to a positive factor, as a canonical
quadratic form of signature $(3,3)$ on the six-dimensional space $\Lambda^{2}
V$. This gives a homomorphism from the identity component of the general linear group $GL^{+}(V)$ to the conformal group of the indefinite form, which
is a local isomorphism. The significance of this is that geometrical structures
on $V$ can be expressed in terms of the six-dimensional space $\Lambda^{2} V$
with its quadratic form. Now if $X$ is an oriented $4$-manifold we can apply
this idea to the cotangent spaces of $X$. Many important differential geometric structures
on $X$ can fruitfully be expressed in terms of the bundle of $2$-forms,  with its quadratic wedge product and exterior derivative. We recall some examples
\begin{itemize}  
\item A {\it conformal structure} on $X$ is given by a 3-dimensional subbundle
 $\Lambda^{+}\subset \Lambda^{2}$ on which the form is strictly positive.
 A
{\it Riemannian metric} is specified by the addition of a choice of volume form.
\item A (compatibly oriented) {\it symplectic structure} on $X$ is a closed
2-form $\omega$ which is positive at every point, $\omega\wedge \omega>0$.
\item An {\it almost-complex structure} on $X$ is given by a $2$-dimensional
oriented
subbundle $\Lambda^{2,0}\subset \Lambda^{2}$ on which the form is strictly
positive.
\item A symplectic form $\omega$ {\it tames} an almost-complex structure
$\Lambda^{2,0}$ if $\langle \omega\rangle + \Lambda^{2,0}$ is a $3$-dimensional
subbundle on which the form is positive, and the induced orientation agrees
with a standard orientation inherited from that of $X$. The symplectic form is {\it compatible}
with the almost-complex structure if $\omega$ is orthogonal to $\Lambda^{2,0}$ and in this case we say that we have an {\it almost-Kahler} structure.
\item A {\it complex-symplectic structure} is given by a pair of closed two
forms $\theta_{1}, \theta_{2}$ such that
$$  \theta_{1}^{2}= \theta_{2}^{2} \  \ ,\ \  \theta_{1}\wedge \theta_{2}=0.
$$
\item A {\it complex structure} is given by an almost-complex structure $\Lambda^{2,0}$
such that in a neighbourhood of each point there is a complex-symplectic
structure whose span is $\Lambda^{2,0}$.
\item A {\it Kahler structure} is an almost-Kahler structure such that $\Lambda^{2,0}$
defines a complex structure.
\item A {\it hyperkahler structure} is given by three closed two-forms $\theta_{1},\theta_{2}, \theta_{3}$ such that
$$ \theta_{1}^{2}= \theta_{2}^{2} = \theta_{3}^{2} \ \ ,\ \ \theta_{i} \wedge
\theta_{j}= 0 \ {\rm for}\ i \neq j. $$
\end{itemize}

A feature running through  these examples is the interaction
between pointwise, algebraic, constraints and the differential constraint
furnished by the exterior derivative. In this article we introduce and begin
the study of a  very general class of questions of this nature, and discuss
the possibility  of various applications to four-dimensional differential geometry.

It is a pleasure to acknowledge the influence of ideas of Gromov and Sullivan
on this article, 
through conversations in the 1980's which in some instances may only have been subconsciously
absorbed by the author at the time. 
 
\section{A class of elliptic PDE}
 Let $P$ be a three-dimensional submanifold in the vector space $\Lambda^{2}
 V$, where $V$ is as in the previous section. We say that
 \begin{itemize}
 \item $P$ has {\it negative tangents} if the wedge product is a strictly negative form
  each tangent space $TP_{\omega}$ for all $\omega$ in $P$.
  \item $P$ has {\it negative chords} if for all pairs $\omega, \omega'\in P$ we have $(\omega-\omega')^{2}\leq 0$ with equality if and only if $\omega=\omega'$.
\end{itemize}
Clearly if $P$ has negative chords it also has negative tangents, but the
converse is not true. Now we turn to our $4$-manifold $X$ and consider a
$7$-dimensional submanifold ${\cal P}$ of the total space of the bundle $\Lambda^{2}$
such that the projection map induces a fibration of ${\cal P}$, so each fibre
is a $3$-dimensional submanifold, as considered before. We say that ${\cal
P}$ has negative tangents, or negative chords, if the fibres do. In any case
we can consider the pair of conditions, for a $2$-form $\omega$ on $X$
\begin{equation} \omega \subset {\cal P} \ \ , \ \ d\omega=0. \label{eq:constraint}\end{equation}
Here, of course, the first condition just means that at each point $x\in X$
 the value $\omega(x)$ is constrained to lie in the given $3$-dimensional
 submanifold ${\cal P}_{x} \subset \Lambda^{2}_{x}$. The pair of conditions
 represent a partial differential equation over $X$. Now suppose that $X$
 is compact and  fix a maximal positive
 subspace $H^{2}_{+}\subset H^{2}(X; {\bf R})$, for example the span of the
 self-dual harmonic forms for some Riemannian metric. Given a class $C\in
 H^{2}(X;{\bf R})$, we augment Equation~\ref{eq:constraint} by the cohomological condition
 \begin{equation}  [\omega] \in C+  H^{2}_{+} \subset H^{2}(X;{\bf R}). \label{eq:cohconstraint}
 \end{equation}

 \begin{prop}
 \begin{itemize}
  \item If ${\cal P}$ has negative chords then there is at most one solution
 $\omega$ of the  constraints (1), (2). 
 \item If ${\cal P}$ has negative tangents and $\omega$ is a solution of
 the constraints (1), (2) then there is a neighbourhood $N$ of $\omega$ (in,
 say, the $C^{\infty}$ topology on $2$-forms) such that there are no other
 solutions of the constraints in $N$. Further, if ${\cal P}^{(t)}$ is a smooth
 $1$-parameter family of deformations of ${\cal P}= {\cal P}^{(0)}$ then
 we can choose $N$ so that
 for small $t$ there is a unique solution of the deformed constraint in $N$.
\end{itemize} 
 \end{prop}
 
 The proof of the first item is just to observe that if $\omega, \omega'$
 are two solutions then $(\omega-\omega')^{2}$ is non-positive, pointwise on $X$,
 by the negative chord assumption. On the other hand the de Rham cohomology
 class of $\omega-\omega'$ lies in $H^{2}_{+}$ so
 $$ \int_{X} (\omega-\omega')^{2} \geq 0. $$
 Thus $(\omega-\omega')^{2}$ vanishes identically and the negative chord
 assumption implies that $\omega=\omega'$.
  
 For the second item we consider, for each point $x$ of $X$ the tangent space
 to the submanifold ${\cal P}_{x}$ at the point $\omega(x)$. This is a maximal
 negative subspace for the wedge product form and so $\omega$ determines a conformal
 structure, on $X$ (the orthogonal complement of the tangent space is a maximal
 positive subspace). For convenience, we fix a Riemannian metric $g$ in this
 conformal class.  The condition that a nearby form $\omega+ \eta$ lies
 in ${\cal P}$ takes the shape
 $$  \eta^{+} = Q(\eta), $$
 where $\eta^{+}$ is the self-dual part of $\eta$ with respect to $g$, and
 $Q$ is a smooth map with $Q(\eta)= O(\eta^{2})$.
 We choose $2$-forms representing the cohomology classes in $H^{2}_{+}$,
 so $H^{2}_{+}$ can be regarded as a finite-dimensional vector space of closed
 $2$-forms on $X$. Then closed forms $\omega+\eta$ satisfying our cohomological
 constraint can be expressed as $\omega+ d a + h$ where $h\in H^{2}_{+}$ and
 where $a$ is a $1$-form satisfying the \lq\lq gauge fixing'' constraint $d^{*}
 a=0$. Thus our constraints correspond to the solutions of the PDE
 $$  d^{*} a=0, d^{+}a =  Q(da+ h) - h^{+}, $$
 where $d^{+}$ denotes the self-dual component of $d$. This is not quite
 a $1-1$ correspondence, we need to identify the solutions $a, a+\alpha$
  where $\alpha$ is a harmonic $1$-form on $X$. 
 The essential point now is that the linear operator
 $$  d^{*} \oplus d^{+}:\Omega^{1} \rightarrow \Omega^{0} \oplus
 \Omega^{2}_{+} $$
 is {\it elliptic}. This means that our problem can be viewed as solving
 a non-linear elliptic PDE and we can apply the implicit function theorem,
 in a standard fashion. The linearisation of the problem is represented by
 the linear map
 $$ L=  d^{*}\oplus d^{+}: \Omega^{1}/{\cal H}^{1} \rightarrow \Omega^{0}/{\cal
 H}^{0} \oplus \Omega^{2}_{+}/ H^{2}_{+}, $$
 where ${\cal H}^{i}$ denotes the space of harmonic $i$-forms, for $i=0,1$.
 We claim that the kernel and cokernel of $L$ both vanish. If $d^{+}a\in H^{2}_{+}$ we have $d d^{+} a=0$ (since the forms in $H^{2}_{+}$
 are closed) and then
 $$ 0= \int_{X} dd^{+} a \wedge a=  \int_{X} d^{+}a\wedge da = \int_{X} \vert d^{+} a \vert^{2}, $$
 so $d^{+}a=0$. Now we use a fundamental identity
 \begin{equation}  0= \int_{X} da \wedge da= \int_{X} \vert d^{+} a \vert^{2} - \vert d^{-}
 a \vert^{2}, \label{eq:identity}\end{equation}
 so $d^{-}a$ vanishes as well. It follows that the kernel of $L$ is trivial.
 The cokernel of $d^{+}: \Omega^{1} \rightarrow \Omega^{+}$ is represented
 by the space of self-dual harmonic forms ${\cal H}^{2}_{+}$. The assertion
 that the cokernel of $L$ is trivial is equivalent to the statement that
 $L^{2}$ projection  $\pi:H^{2}_{+}\rightarrow {\cal H}^{2}_{+}$ is surjective. Since these
 are both maximal positive subspaces for the cup product form they have the
 same dimension, so it is equivalent
 to prove that $\pi$ is injective. But if a form $h\in H^{2}_{+}$ is in the kernel of
 $\pi$ the cohomology class $[h]$ lies in the subspace ${\cal H}^{2}_{-}\subset
 H^{2}(X;{\bf R})$ defined by the anti-self dual harmonic forms, so $\int_{X}h^{2}\leq
 0$ and it follows that $h=0$.
 
 To sum up, for a given choice of ${\cal P}$ our constraints are represented
 by a system of nonlinear elliptic PDE whose linearisation is invertible. Now the
 assertions in the second item follow in a standard way from the implicit
 function theorem.
 
 \
 For the rest of this paper we will consider constraint manifolds ${\cal
 P}$ with negative tangents.
  Notice that we could generalise the set-up slightly by choosing a submanifold
 $Q\subset H^{2}(X;{\bf R})$ with the property that at each point 
 the tangent space of $Q$ is maximal positive subspace for the cup product
 form. Then we can take in place of the constraint (2) the condition
 $  [\omega]\in Q$. The proof above goes through without essential change.
 We could also express things differently by considering the moduli space
 ${\cal M}$ of solutions to (1) with no cohomological constraint. Then the
 same proof shows that ${\cal M}$ is a manifold of dimension $b^{2}_{-}(X)$
 and that the map $\omega \mapsto [\omega]$ defines an immersion of ${\cal
 M}$ in $H^{2}(X;{\bf R})$ whose derivative at each point takes the tangent
 space of ${\cal M}$ to   a maximal
 {\it negative} subspace for the intersection form.

 \section{Examples}
 
 1. We fix a Riemannian metric on $X$ and let ${\cal P}$ be the vector subbundle
 of anti-self dual forms.  Then the solutions of our constraints are just
 the anti-self dual harmonic forms. This case is not very novel but, as we
 have seen in the 
  previous section, is the model for the general situation  for the purposes of   local deformation theory.
 
 2. Take $X={\bf R}^{4}= {\bf R}^{3} \times {\bf R}$ with co-ordinates $(x_{1}, x_{2},
 x_{3}, t)$ and consider a case where ${\cal
 P}$ is preserved by translations in the four variables, so is determined
 by a single $3$-manifold $P\subset \Lambda^{2}({\bf R}^{4})$. (Of course $X$
 is not compact here, but we only want to illustrate the local PDE aspects.)
 In the usual way, we write a $2$-form as a pair of vector fields  $(E,B)$
 $$ \omega = \sum E_{i} dx_{i} \wedge dt + \frac{1}{2}\sum\epsilon_{ijk} B_{i} dx_{j} \wedge dx_{k};
 $$
 that is, we are identifying $\Lambda^{2} {\bf R}^{4}$ with ${\bf R}^{3} \oplus
 {\bf R}^{3}$. The condition that a submanifold $P\subset \Lambda^{2} {\bf
 R}^{4}$ has negative tangents implies that its tangent space at each point
 is transverse to the two ${\bf R}^{3}$ factors, so locally $P$ can be written
 as the graph of a smooth map $F:{\bf R}^{3} \rightarrow {\bf R}^{3}$. That
 is, locally around a given solution, we can write the constraint as $B=F(E)$.
 Consider solutions which are independent of translation in the $t$ variable.
 The condition that $\omega$ is closed  becomes
 $$   \nabla . B=0 , \ \nabla \times E=0. $$
  We can write $E=\nabla u$ for a function $u$ on ${\bf R}^{3}$, so our constraint
  is a nonlinear elliptic PDE for a function $u$ on ${\bf R}^{3}$ of the form
  $$    \nabla.( F(\nabla u))=0. $$
 Let
 $$ ( H_{ij})= \left( \frac{\partial F_{i}}{\partial \xi_{j}}\right), $$ be
 the matrix of derivatives of $F$. The condition that $P$ has negative tangents
 becomes $H+H^{T}>0$ and the linearisation of the problem is the linear elliptic
 PDE
 $$   \sum_{ij} \frac{\partial}{\partial x_{i}} \left(  H_{ij} \frac{\partial
 f}{\partial x_{j}}\right) = 0. $$
 
 \
 
 3.
 The next example is the central one in this article. We fix a volume form
  on our $4$-dimensional real vector space $V$ and also a complex structure,
 i.e. a $2$-dimensional positive subspace $ \Lambda^{2,0}\subset \Lambda^{2} V$. We define $P$ to be the set of positive $(1,1)$-forms whose square is the given volume
form. Then $P$ is one connected component of the set 
$$   \{ \omega \in \Lambda^{2} V: \omega^{2}=1, \omega \wedge \Lambda^{2,0}=0\},
$$
(the component being fixed by orientation requirements).  
It is easy to see that $P$ has negative chords. In fact if $\omega,\omega'$
are two points in $P$ we can choose complex coordinates $z_{1}, z_{2}$ so
that
$$ \omega = i dz_{1} \wedge d\overline{z}_{1} + i dz_{2} \wedge d\overline{z}_{2}\
\ ,\ \  \omega'= \lambda i dz_{1} \wedge d\overline{z}_{1} + \lambda^{-1} i
dz_{2} \wedge d\overline{z}_{2}, $$
where  $\lambda>0$. Then $$(\omega-\omega')^{2}= (\lambda-1)(\lambda^{-1}-1)=
2-(\lambda+\lambda^{-1}) \leq 0,$$
with equality if and only if $\lambda=1$. 

Thus if $X$ is a $4$-manifold with a volume form and a choice of almost-complex
structure we get a constraint manifold ${\cal P}$ with negative chords, and
a problem of the kind we have been considering. In particular, we can consider
the case when the almost-complex structure is integrable. Then our problem
becomes the renowned Calabi conjecture (in the case of two complex variables)
solved by Yau: prescribing the volume form of a Kahler metric. The solution
can be most easily expressed in terms of the moduli space ${\cal M}$ of solutions
to (1): it maps bijectively to the intersection of the  Kahler cone in
$ H^{1,1}(X)$ with the quadric $[\omega]^{2}={\rm Vol}$, where Vol
 is the
integral of the prescribed volume form. Notice however that the problem formulated
in terms of the cohomological constraints (2) does {\it not} always have a
solution. In the most extreme case, we could take $X$ to be a complex surface
which does not admit any Kahler metric. More generally, by deforming our
choice of $C$ and $H^{2}_{+}$ we can deform from a case when a solution exists
to a case when it does not, and understanding this phenomenon is essentially
the question of understanding the boundary of the Kahler cone. 
 
 The extension of the Calabi-Yau theory to the case when the almost structure
 is not integrable has been considered recently by Weinkove \cite{kn:W}. Suppose,
 for simplicity, that $b^{2}_{+}(X)=1$ and suppose that $\omega_{0}$ is a symplectic
 form compatible with the given almost-complex structure. In this case we
 take $C=0$ and  $H^{2}_{+}$ to be the $1$-dimensional space spanned by $\omega_{0}$.
 Then the cohomology class of any solution of (1), (2) is fixed by the prescribed
 volume form and without loss of generality we can suppose it is the same
 as $[\omega_{0}]$. Weinkove extends Yau's {\it a priori} estimates to prove existence
 under the assumption that the Nijenhius tensor of the almost complex structure
 is small in a suitable sense. 
 
 In general, we will say that a  constraint ${\cal P}\subset \Lambda^{2}X$ is {\it unimodular} if there is a  volume form $\rho$ on $X$ such that $\omega^{2}=\rho$
for any section $\omega \subset {\cal P}$. So for a unimodular constraint
any solution of (1) is a symplectic form.
 
  \

 4. Our final example is in some ways a simple modification of the previous one. We take an almost-complex structure  and a  fixed positive $(1,1)$-form $\Theta$ and we consider the submanifold ${\cal P}$ of positive $(1,1)$ forms
$\omega$ with $(\omega-\Theta)^{2}$ a given volume form on $X$. In each fibre
this is a translate of the submanifold considered before, so again has negative
chords. 
This example arises in the following way. Let $X$ be a hyperkahler $4$-manifold,
with an orthogonal triple of closed forms $\theta_{1}, \theta_{2}, \theta_{3}$
and let $\omega$ be another symplectic form on $X$. We define three functions
$\mu_{i}=\mu_{i}(\omega)$ on $X$ by
$$  \mu_{i} = \frac{\theta_{i} \wedge \omega}{\omega^{2}}. $$
These arose in \cite{kn:D} as the \lq\lq moment maps'' for the action of the symplectomorphism
group of $(X,\omega)$ with respect to the $\theta_{i}$ and the triple $(\mu_{1},
\mu_{2}, \mu_{3})$ can be regarded as a hyperkahler moment map. Now we ask
the question: given three functions $f_{i}$ on $X$, can we find a symplectic
form $\omega$ with $\mu_{i}(\omega)=f_{i}$? An obvious necessary condition
is that $F=\sqrt{\sum f_{i}^{2}}$ does not vanish anywhere on $X$ (for then the
self-dual part of $\omega$ would vanish and $\omega^{2}$ would be negative). Assuming this condition, we write $\sum f_{i} \theta_{i} = F \sigma$ where
$\sigma$ is a unit self-dual $2$-form. Then $\sigma$ determines an almost-complex structure
on $X$ ( with $\Lambda^{2,0}$ equal to the orthogonal complement of $\sigma$
in the span of $\theta_{1}, \theta_{2}, \theta_{3}$) and $\sigma$ is a positive
form  of type
(1,1) with respect to this structure.  It is easy to check that the condition
that $\mu_{i}(\omega)=f_{i}$ can be expressed in the form above, with
$\Theta= \sigma/2F$ and the volume form $\sigma^{2}/(4 F^{2})$.

\section{Partial regularity theory}

We have introduced a very general class of elliptic PDE problems on $4$-manifolds
and shown that their behaviour with respect to deformations is straightforward.
The crux of the matter, as far as proving existence results goes, is thus
to obtain {\it a priori} estimates for solutions. Here we make some small steps in this direction, assuming an $L^{\infty}$ bound. So throughout this section we
assume that ${\cal P}$ is a constraint manifold with negative tangents and
$\omega$ is a closed form in  ${\cal P}$ with $\vert \omega \vert \leq K$ at
each point. Here the norm is defined by some fixed auxiliary Riemannian metric
$g_{0}$ on $X$. We will also use the conformal structure determined by the
pair $\omega, {\cal P}$ which for convenience we promote to another metric $g$,
with, say, the same volume form as $g_{0}$. Since the set
$ \{ \Omega\in {\cal P}, \vert \Omega \vert \leq K\}$  is compact the metrics
$g,g_{0}$ are uniformly equivalent, for fixed $K$.  
 
 Our estimates will depend on $K$ and ${\cal P}$. More precisely,
the dependence on ${\cal P}$ will involve local quantities that could be
written down explicitly and, crucially, will be uniform with respect to continuous
families of constraints ${\cal P}_{t}$ (with fixed $K$).  

\begin{lem} There exists $C=C(K,{\cal P})$ such that
$$ \Vert \nabla \omega \Vert_{L^{2}} \leq C $$
 \end{lem}
 (Here $\nabla$ is the covariant derivative associated with the fixed metric
 $g_{0}$.)
 
 It suffices to show that for any vector field $v$ on $X$ the Lie derivative
 $L_{v} \omega$ is bounded in $L^{2}$ by some $C(K,{\cal P}, v)$. For then
 we consider a cover of $X$ by $D$ coordinate patches and $4D$ vector fields
 such that on each patch four of the vector fields are the standard constant
 unit fields, for which the Lie derivative becomes the ordinary derivative.
 Write $\omega_{v}$ for $L_{v} \omega$. So $d\omega_{v}=0$, and in fact $\omega_{v}=
 d i_{v}(\omega)$ is exact.
 Imagine first that the flow generated by $v$ preserves the constraint manifold
 ${\cal P}$. Then
 if we apply the Lie derivative to the condition $\omega \subset {\cal P}$
 we obtain an identity 
 $$  (  \omega_{v})^{+}=0, $$
 where $()^{+}$ denotes the self-dual part with respect to metric $g$
  defined by $\omega$ and ${\cal P}$. In general we will have an
 identity $$ (\omega_{v})^{+}= \rho, $$
 where the $L^{\infty}$ norm of $\rho$ can be bounded in terms of of $v,
 {\cal P}$ and  $K$. In particular, the $L^{2}$ norm of
$(\omega_{v})^{+}$ is {\it a priori} bounded (since the metrics $g$ and
$g_{0}$ are uniformly equivalent we can take the $L^{2}$ norm with respect
to either here). Now using the metric $g$ we have, since $\omega_{v}$ is exact
$$ \Vert \omega_{v} \Vert^{2}_{L^{2}} = 2 \Vert (\omega_{v})^{+}\Vert^{2}_{L^{2}}
$$
just as in (3), and we obtain the desired $L^{2}$ bound on $\omega_{v}$.

\begin{lem}
There is a constant $c$, depending on $K$ and ${\cal P}$ such that if 
$$ \int_{B(r)} \vert \nabla \omega \vert^{2} \leq c r^{2}, $$
for all $r$-balls $B(r)$ in $X$ with $r\leq r_{0}$ then for any $p$
$$ \Vert \nabla \omega \Vert_{L^{p}} \leq C, $$
where $C$ depends on $p,K, {\cal P}, r_{0}$.
\end{lem}

We use of some of the ideas developed in \cite{kn:DS}. Fix attention on balls
with a given centre,  and choose local coordinates about this point such that the metric $g$ is the standard Euclidean metric at the origin.
Let $\Lambda^{\pm}_{0}$ be the space of $\pm$ self-dual forms for the euclidean
structure.
The metric $g$ is represented by a tensor $\mu\in Hom(\Lambda^{-}_{0}, \Lambda^{+}_{0}$, vanishing at the origin,
such that the $g$-anti-self-dual $2$-forms have the shape $\sigma+\mu(\sigma)$ for
$\sigma$ in $\Lambda^{-}_{0}$. Now for any such tensor field $\mu$, defined over
the unit ball $B$  in ${\bf R}^{4}$ say, consider the operator
$$ d^{+}_{\mu}=  d^{+}_{0} + \mu d^{-}_{0}: \Omega^{1} \rightarrow \Omega^{+}. $$
Here $d^{+}_{0}, d^{-}_{0}$ are the constant co-efficient operators defined by $\Lambda^{+}_{0},
\Lambda^{-}_{0}$. 
The basic point is that, for any given $p$, if $\mu$ is sufficiently small
in $L^{\infty}$
 the usual Calderon-Zygmund theory can be applied
to this operator, regarded as a perturbation of $d^{+}_{0}$. More precisely,
for given $p$ there is a $\delta=\delta(p)$ such that if $\vert \mu \vert\leq
\delta$ then for any 
 $1$-form $\alpha$ over the ball with $d^{*}\alpha=0$ we have
$$ \Vert d \alpha \Vert_{L^{p}(B/2)} \leq C_{p} (\Vert d^{+}_{\mu} \alpha\Vert_{L^{p}(B)}+
\Vert d\alpha \Vert_{L^{2}(B)}). $$
It follows that for closed $2$-forms  $\rho$  over $B$ we have an inequality
\begin{equation} \Vert \rho\Vert_{L^{p}(B/2)} \leq C_{p} (\Vert \rho^{+,\mu} \Vert_{L^{p}(B)}
+ \Vert \rho \Vert_{L^{2}(B)})\label{eq:CZ} \end{equation} 
where $\rho^{+,\mu}$ denotes the self-dual part with respect to the conformal
structure defined by $\mu$. (To see this we write $\rho=d\alpha$ with $d^{*}
\alpha=0$. )
To apply this idea in our situation we first fix some $p>4$.
For each $r$ we rescale the ball $B(r)$ to the unit ball in ${\bf R}^{4}$
and
apply (4) to the ordinary derivatives of the form corresponding to $\omega$, defined with respect to
the Euclidean coordinates, which are closed $2$-forms. Since the tensor $\mu$
which arises is determined by the tangent space of ${\cal P}$, which varies
continuously with $\omega$, there is some $\epsilon
= \epsilon(p)$ such that (4) holds provided the oscillation of $\omega$ over
$B(r)$ is less than $\epsilon$. (Here the \lq\lq oscillation''of $\omega$
refers to the oscillation of the coefficients in the fixed coordinate system.)
 Transferring back to $B(r)$, and taking account of the
rescaling behaviour of the quantities involved,  we obtain an inequality of
the form
$$  r^{-4/p} \Vert \nabla \omega \Vert_{L^{p}(B(r/2))} \leq C_{p} 
r^{-2} \Vert \nabla \omega \Vert_{L^{2}(B(r))} +  C'_{p} , $$
for constants $C_{p},C'_{p}$ depending on $p, K, {\cal P}$. This holds provided
the oscillation of $\omega$ over $B(r)$ is less than $\epsilon$. On the other
hand the Sobolev inequalities tell us that the oscillation of $\omega$ over
$B(r/2)$ is bounded by  a multiple of $r^{(p-4)/p}\Vert \nabla \omega \Vert_{L^{p}(B(r/2))}$. So we conclude that
if the oscillation of $\omega$ over $B(r)$ is less than $\epsilon$ then the
oscillation of $\omega$ over $B(r/2)$ is at most
$$   C \left( r^{-2} \int_{B(r)} \vert\nabla \omega\vert^{2} \right)^{1/2} + C' r \leq C c^{1/2} + C' r. $$
Thus if $c$ is sufficiently small we can arrange that the oscillation of
$\omega$ over the half-sized ball is less than $\epsilon/10$, say, once $r$ is small
enough. Applying this to a pair of balls, we see that there is some $r_{0}$
such that if the oscillation of $\omega$ over
all $r$-balls is less than $\epsilon$, for all $r\leq r_{0}$, then this oscillation
is actually less than $\epsilon/2$. It follows then, by a continuity argument
taking $r\rightarrow 0$, that the oscillation can never exceed $\epsilon$
over balls of a fixed small size and this gives an  {\it a priori} bound on the $L^{p}$ norm of $\nabla \omega$ for this fixed $p$. Now we have a fixed
Holder bound on $\omega$ and we can repeat the discussion, starting with this,
to get an $L^{p}$ bound on $\nabla \omega$ for any $p$.

\

\

Let $B(r)$ be an embedded geodesic ball in $(X,g_{0})$. We let $\hat{\omega}$
denote the $2$-form over $B(r)$ obtained by evaluating $\omega$ at the centre
of the ball and
extending over the ball by radial parallel transport along geodesics.
\begin{lem} For any $c>0$
there is a constant $\gamma$, depending on $K,{\cal P},c $ such that if for
all $\rho$-balls $B(\rho)$ (with $\rho<\rho_{0}$ ) we have
$$  \int_{B(\rho)} \vert \omega - \hat{\omega} \vert^{2} \leq \gamma \rho^{4} $$
then the hypothesis of Lemma 2 is satisfied, i.e. $$\int_{B(r)}\vert \nabla
\omega\vert^{2} \leq c r^{2}$$ for all $r<r_{0}(\rho_{0},K, {\cal P})$.
\end{lem}

As before, we work in small balls centred on a given point  in $X$ and standard
 coordinates about this point. There is no loss
in supposing that the metric $g_{0}$ is the Euclidean metric in these coordinates:
then the form $\hat{\omega}$ is just obtained
by freezing the coefficients at the centre of the ball. Let $v$ be a 
constant vector field in these Euclidean co-ordinates. Fix a cut-off function
$\beta$ equal to $1$ on $[0,1/2]$ and supported in $[0,1)$ and let $\beta_{r}$
be the function on $X$ given by $\beta_{r}(x)= \beta(\vert x\vert/r)$ in
these Euclidean coordinates.  Set 
$$  I(r)= \int_{B(r)} \beta_{r}  \omega_{v}\wedge \omega_{v} . $$
As in the proof of Lemma 1, it suffices to show that, when  $r$ is small,  $I(r)\leq c r^{2}$ (for a
different constant $c$). Suppose first that $r=1$ and write $\omega= \hat{\omega}
+\eta$, so $\Vert \eta \Vert^{2}_{L^{2}(B(1))} \leq \gamma$. By an $L^{2}$ version
of the Poincar\'e inequality, as in \cite{kn:U}, we can find a $1$-form $\alpha$
over the ball with $\eta = d\alpha $ and $ \Vert \nabla \alpha \Vert^{2}_{L^{2}(B(1))} \leq C \gamma$
for some fixed constant $C$. Now we have $\omega_{v}= \eta_{v}$, since the
form $\hat{\omega}$ is constant in these coordinates, and $\omega_{v} = d
\alpha_{v}$, where $\Vert\alpha_{v} \Vert^{2}_{L^{2}(B(1))} \leq C \gamma$.
Then 
$$  I(1)= \int_{B(1)} \beta_{1} \omega_{v} \wedge d \alpha_{v}= \int_{B(1)}
d\beta_{1}\wedge \omega_{v} \wedge \alpha_{v}, $$ so
$$ I(1) \leq \Vert d \beta_{1} \Vert_{L^{\infty}} \Vert \omega_{v} \Vert_{L^{2}(B(1))}
\Vert \alpha_{v} \Vert_{L^{2}(B(1))} . $$
Using the $L^{2}$ bound on $\alpha_{v}$ and the fact that the $L^{2}$ norm
of $\omega_{v}$ on $B(1)$ is controlled by $\sqrt{I(2)}$ we obtain
$$   I(1) \leq C \sqrt{\delta} \sqrt{I(2)}. $$
For general $r$, we scale the $ 2r$ ball to the unit ball and apply the same
argument. Taking account of the scaling behaviour of the various quantities
involved one obtains an inequality
$$  I(r) \leq (C \sqrt{\gamma}) r \sqrt{I(2r)}. $$
It is elementary to show that this implies the desired decay condition on
$I(r)$ as $r$ tends to $0$. If we put  $J(r) = r^{-2} I(r)$ then $J(r)\leq
2 C \sqrt{\gamma} \sqrt{ J(2r)}$. So if $L_{n} = \log J(2^{-n})$ we have
$$ L_{n+1} \leq  \frac{L_{n}}{2} + \sigma$$
where $\sigma = \log(2C \sqrt{\gamma})$.  This gives
$$ L_{n} \leq 2^{-n} (L_{0}- \sigma) + 2 \sigma. $$
Combined with the global {\it a priori} bound of Lemma 1, which controls
$L_{0}$,  this yields the
desired result. 

To sum up we have
\begin{prop}
Suppose ${\cal P}$ is a constraint manifold with negative tangents and $K>0$.
There is a $\gamma=\gamma({\cal P}, K)$ with the following property. For each $k,p,
r_{0}$
there is a constant $C=C(k,K, {\cal P}, p,r_{0})$ such that if  $\omega$ is
any solution of (1)
with $\Vert \omega\Vert_{L^{\infty}} \leq K$
and with $$ \int_{B(r)} \vert \omega-\hat{\omega}\vert^{2}\leq \gamma r^{4}$$ for
all $r\leq r_{0}$  then $\Vert \omega \Vert_{L^{p}_{k}} \leq C. $
\end{prop}

For $k=1$ this is a combination of the two preceding Lemmas. The extension
to higher derivatives follows from a straightforward bootstrapping argument.

\

\

Of course the hypotheses of Proposition 2 are satisfied if $\omega$ is bounded
in $L^{\infty}$ and has any fixed modulus of continuity, such as a Holder
bound. 

\section{Discussion}
  \subsection{Motivation from symplectic topology}
  We will outline the, rather speculative, possibility of applications of
  these ideas to questions in symplectic topology. Recall that the fundamental
  topological invariants of a  symplectic form $\omega$ on a compact $4$-manifold
  $X$ are the first Chern class $c_{1}\in H^{2}(X;{\bf Z})$ and the de Rham class $[\omega]\in
  H^{2}(X;{\bf R})$.
  \begin{question}
  Suppose $X$ is a compact Kahler surface with Kahler form $\omega_{0}$. If
  $\omega$ is any other symplectic form on $X$, with the same Chern class
 and with $[\omega]=[\omega_{0}]$, is there a diffeomorphism $f$ of $X$ with
 $f^{*}(\omega)=\omega_{0}$ ?
 \end{question}
 (McMullen and Taubes \cite{kn:McT} have given examples of inequivalent symplectic structures
  on the same differentiable $4$-manifold, but their examples have different
  Chern classes.)
  
  A line of attack on this  could run as follows. Suppose first that $b^{2}_{+}(X)=1$.
  Given a general symplectic form $\omega$ we choose a unimodular constraint manifold
  ${\cal P}_{1} $ containing it and deform ${\cal P}_{1}$ through a 1-parameter
  family ${\cal P}_{t}$ for $t\in [0,1]$ to a standard Calabi-Yau constraint
  ${\cal P}_{0}$ containing $\omega_{0}$. We suppose that ${\cal P}_{t}$
  are all unimodular, with the same fixed volume form. Then we choose the
  cohomological constraint by taking $C=0$ and  $H^{2}_{+} = \langle \omega \rangle$. Thus for any $t$ the cohomology class of a solution is forced to
be the fixed class $[\omega_{0}]$. At $t=0$ we know that the solution $\omega_{0}$
is unique. If we imagine that we have obtained suitable {\it a priori} estimates
for solutions in the whole family of problems it would follow from the deformation
result, Proposition 1,  that there is a unique solution $\omega_{t}$ for each $t$, varying
smoothly with $t$, and $\omega=\omega_{1}$.  Then Moser's theorem would imply that $\omega$ and $\omega_{0}$ are equivalent forms.

This strategy can be extended to the situation where $b^{2}_{+}>1$. We take $H^{2}_{+}= \langle \omega_{0}\rangle + H^{2,0}$ where
$H^{2,0}$ consists of the real parts of the holomorphic $2$-forms. There
are two cases. If the Chern class $c_{1}$ vanishes in $H^{2}(X;{\bf R})$
 then for any
  $\theta\in H^{2,0}$ and $s\in {\bf R}$ the form $\theta+ s \omega_{0}$ is  symplectic, provided $\theta$ and $s$ are not both zero . (For the zero
set of $\theta$ is either empty or a nontrival  complex curve in $X$, and
the second possibility is excluded by the assumption on $c_{1}$.)  If $c_{1}$ does
not vanish
then the form is symplectic provided that $s\neq 0$. We follow the same procedure
as before and (assuming the {\it a priori} estimates) construct a path $\omega_{t}$
from $\omega_{0}$ to $\omega=\omega_{1}$ with  $[\omega_{t}]= s_{t} [\omega_{0}]
+ [\theta_{t}]$ for $\theta_{t} \in H^{2,0}$. Now $[\omega_{t}]$ is not identically
zero, so in the case $c_{1}=0$ we have a family of \lq\lq standard'' symplectic
forms
$$ \tilde{\omega_{t}} = s_{t} \omega_{0} + \theta_{t}, $$
with $[\tilde{\omega}_{t}]=[\omega_{t}$. Then a version of Moser's Theorem
yields a family of diffeomorphisms $f_{t}$
 with $f_{t}^{*}(\tilde{\omega}_{t})=\omega_{t}$. Since, by hypothesis,
 $\tilde{\omega}_{1}=\omega_{0}$ and $\omega_{1}=\omega$,  the diffeomorphism $f_{1}$ solves our problem. When $c_{1}\neq 0$ the same argument works provided
we know that $s_{t}$ is never $0$. But here we can use one of  the deep results of
Taubes \cite{kn:T}. If $s_{t}$ were zero then there would be a class $\theta=\theta_{t}$
in $H^{2,0}$ which is the class of a symplectic form with the same first
Chern class $c_{1}(X)$. But Taubes' inequality would give
$  c_{1} . [\theta]<0$ whereas in our case $c_{1}. [\theta]=0$ since
$c_{1}$ is represented by a form of type (1,1) and $\theta$ has type $(2,0)+(0,2)$. 

\subsection{The almost-complex case.}
We have seen that, even within the standard framework of the Calabi-Yau problem
on complex surfaces, solutions can blow up. In that case, everything can
be understood in terms of the class $[\omega]$ and the Kahler cone. These
difficulties become more acute in the nonintegrable situation. If $J_{0}$
is any almost complex structure on a $4$-manifold $X$ we  can construct
a $1$-parameter family  $J_{t}$ such that there is no symplectic form, in
any cohomology class, compatible
with $J_{1}$. (This is in contrast with the integrable case, where deformations
of a Kahler surface are Kahler). To do this we can simply deform $J_{0}$ in a small neighbourhood
of a point so that $J_{1}$ admits a null homologous pseudo-holomorphic curve.
More generally we could consider null-homologous  currents $T$ whose (1,1) part
is positive. Thus if we form a $1$-parameter family of constraints ${\cal P}_{t}$
using these $J_{t}$, and any fixed volume form, solutions $\omega_{t}$ must
blow up sometime before $t=1$, however we constrain   the cohomology class
$[\omega_{t}]$. It seems plausible that solutions blow-up at the first time
$t$ when  a null-homologous 
 (1,1)-current $T$ appears, and become singular along the support of
 $T$. It is a result of Sullivan \cite{kn:Sul} that the nonexistence
of such currents implies   the existence of a symplectic
form {\it taming} the almost-complex structure. These considerations lead
us to  formulate a tentative
conjecture.
\begin{conj}
Let $X$ be a compact $4$-manifold and let $\Omega$
be a symplectic form on $X$.
If ${\cal P}$ is a constraint manifold defined by an almost-complex
 structure which is tamed by $\Omega$, and any smooth volume form, then
 there are $C^{\infty}$ {\it a priori} bounds on a closed form $\omega \subset
 {\cal P}$ with $[\omega]=[\Omega]$.
 \end{conj}
(By the results of Weinkove \cite{kn:W} it suffices to obtain $L^{\infty}$
bounds on $\omega$.)
 
 This conjecture  is relevant to the following problem.
 \begin{question}
 If $J$ is an almost-complex structure on a compact $4$-manifold which is
 tamed by a symplectic form, is there a symplectic form compatible with $J$?
 \end{question}
 
 If Conjecture 1 were true it would imply, by a simple deformation argument,an
 affirmative answer to Question 2 in the case when $b^{2}_{+}=1$, see the
 discussion in \cite{kn:W}. Such a result
 would be of interest even in the integrable case. It is a well-known fact
 that any compact complex surface with $b^{2}_{+}$ odd is Kahler. This was originally
 obtained from the classification theory, but more recently direct proofs
 have been given \cite{kn:B}, \cite{kn:L}. Harvey and Lawson showed that
 any surface with $b^{2}_{+}$ odd admits a symplectic form taming the complex
 structure (\cite{kn:HL}, Theorem 26 and page 185), so  a positive answer
 to the question above would yield another  proof of
 the Kahler property, in the case when $b^{2}_{+}=1$. (Turning things around,
 one might hope that the techniques used in \cite{kn:B}, \cite{kn:L} could have
 some bearing on Conjecture 1).
 
 In connection with Question 2, note first that this is special to 4-dimensions.
 In  higher dimensions a generic almost-complex structure does not admit
 any compatible symplectic stucture, even locally. In another direction,
 the answer is known to be affirmative in the case when the taming form is
 the standard symplectic form on ${\bf C P}^{2}$, by an argument of Gromov.
 This
  constructs a compatible form by averaging over the currents furnished
 by pseudoholomorphic spheres.
 
 \subsection{Hyperkahler structures}
 Recall that complex-symplectic and hyperkahler structures on $4$-manifolds
 can be described by, respectively, pairs and triples of orthonormal closed
 $2$-forms. We can ask
 \begin{question} Let $X$ be a compact oriented $4$-manifold
 \begin{itemize}
 \item Suppose there are closed two-forms $\theta_{1}, \theta_{2}$ on $X$
 such that each point $\theta_{i}$ span a positive $2$-plane in $\Lambda^{2}$.
 Does $X$ admit a complex-symplectic structure?
\item Suppose there are closed two-forms $\theta_{1}, \theta_{2}, \theta_{3}$
on $X$ such that at each point $\theta_{i}$ span a positive $3$-plane in
$\Lambda^{2}$. Does $X$ admit a hyperkahler structure?
 \end{itemize}\end{question}
 The hypotheses are equivalent to saying that the symmetric matrices
 ($2\times 2$ and $3\times 3$ in the two cases)
  $\theta_{i}\wedge
 \theta_{j}$ are positive definite, and the question asks whether we can
 find another choice of closed forms $\tilde{\theta_{i}}$ to make the corresponding
 matrix $\tilde{\theta_{i}}\wedge \tilde{\theta_{j}}$ the identity.
 
 For simplicity we will just discuss the second version of the question,
 for triples of forms and in the case when $X$ is simply connected. The
  hypotheses imply that the symplectic structure $\theta_{1}$, say, has first
  Chern class zero, and a result of Morgan and Szabo \cite{kn:MS} tells us
   that $b^{2}_{+}(X)=3$. Thus the $\theta_{i}$ generate a maximal positive
  subspace $H^{2}_{+}$. For reasons that will appear soon we  make
  an additional auxiliary assumption. We suppose that there is an
   involution $\sigma:X\rightarrow X$ with $$\sigma^{*}(\theta_{1})=\theta_{1}\
   ,\
   \sigma^{*}(\theta_{2})=-\theta_{2}\ ,\ \sigma^{*}(\theta_{3})=-\theta_{3}.
   $$ 
   (A model case of this set-up is given by taking $X$ to be a K3 surface
    double-covering  the plane, with $\sigma$ the covering involution.)
   We show that under this extra assumption the truth of Conjecture 1 would
   imply a positive answer to our question. To see this we let $U_{0}$ be the
   orthogonal complement  of $\theta_{1}$ in the span of the $\theta_{i}$,
   so $U_{0}$ defines an almost-complex structure on $X$, compatible with
   $\theta_{1}$. We can chooose a smooth
   homotopy $U_{t}$ to the subbundle $U_{1}= {\rm Span}(\theta_{2}, \theta_{3})$
  giving a $1$-parameter family of almost-complex structures, all tamed by
  $\theta_{1}$. We choose a $1$-parameter family of volume forms equal to
  $\omega_{1}^{2}$ when $t=0$ and to $\omega_{2}^{2}$ when $t=1$. We make
  all these choices $\sigma$-invariant. Thus we get a $1$-parameter family
  of constraints ${\cal P}_{t}$ and seek $\omega_{t}$ with $[\omega_{t}]$
  in $H^{2}_{+}$. The form $\theta_{1}$ gives a solution at time $0$, so
  $\omega_{0}=\theta_{1}$. The $\sigma$-invariance of the problem means that
  in this case we can fix the cohomology class of $\omega_{t}$ to be a multiple
  of the taming form $\theta_{1}$ and still solve the local deformation problem.
   Thus the hypotheses
 of Conjecture 1 are fulfilled so, assuming the conjecture, we can continue
 the solution over the whole interval to obtain an $\omega_{1}$ which satisfies
$$\omega_{1}^{2}=\theta_{2}^{2}\ \ ,\ \ \omega_{1}\wedge \theta_{2}=\omega_{1}\wedge
\theta_{3}=0. $$ 
  Now $\omega_{1}, \theta_{2}$ define a complex-symplectic stucture on $X$,
  and it well known that under our hypotheses $X$ must be a K3 surface and  hence
  hyperkahler.

  The only purpose of the involution $\sigma$ in this argument is to fix
  the cohomology class of the form in the $1$-parameter family. Of course
  one could hope that a suitable extension of the ideas could remove this
  assumption, but our discussion is only  intended to illustrate how these
  ideas might possibly be useful.
  
  \subsection{Relation with Gromov's Theory}
  We have seen that there is a \lq\lq Calabi-Yau'' constraint manifold 
  ${\cal P}$ associated to an almost complex-structure $J$ (and volume form)
  on a $4$-manifold. On the other hand, we have the notion of a \lq\lq J-holomorphic''
  curve, with renowned applications in global symplectic geometry due to
  Gromov and others. These ideas can be related, and extended,   as we will
  now explain. Consider the Grassmannian $Gr_{2}(\bf R^{4})$ of oriented $2$-planes.
  It can be identified with the space of null rays in $\Lambda^{2} {\bf R}^{4}$
  by the map which takes a rank two $2$-form to its kernel. Thus, choosing
  a Euclidean metric on ${\bf R}^{4}$ and hence a decomposition into self-dual
  and anti-self dual forms, the Grassmannian is identified with pairs
  $(\omega_{+}, \omega_{-})$ where $\vert \omega_{+}\vert^{2} = \vert \omega_{-}\vert^{2}=1$,
  or in other words with the product $S^{2}_{+}\times S^{2}_{-}$ of the unit
  spheres in the $3$-dimensional vector spaces $\Lambda^{2}_{\pm}$.  Now
  there is a canonical conformal structure of signature $(2,2)$ on $Gr_{2}({\bf
  R}^{4})$, induced by
  the wedge product form and in  \cite{kn:Gr}  Gromov considers embedded surfaces
  $T\subset  Gr_{2}({\bf R}^{4})$ on which the conformal structure is negative
  definite. A prototype for such a surface is given by  $ T(\omega_{+})= \{ \omega_{+}\} \times
  S^{2}_{-}$ for some fixed $\omega_{+}$ in $S^{2}_{+}$. This is just the
   set of $2$-planes which are complex subspaces with respect to the almost-complex structure defined by $\omega_{+}$ (with $\Lambda^{2,0}$ the orthogonal complement
of $\omega_{+}$ in $\Lambda^{2}_{+})$. On the other hand, Gromov shows that {\it any}
surface $T$ leads to an elliptic equation generalising that defining holomorphic
curves in ${\bf C}^{2}$. Slightly more generally still, let $X$ be an oriented
$4$-manifold and form the bundle $Gr_{2}(X)$ of Grassmannians of $2$-planes
in the tangent spaces of $X$.  Let ${\cal T}$ be a $6$-dimensional submanifold
of the Grassman bundle, fibering over $X$ with each fibre ${\cal T}_{x}$
a \lq\lq negative'' surface in the above sense. Then we have the notion of
a \lq\lq $T$-pseudoholomorphic curve'' in $X$: an immersed surface whose
tangent spaces lie in ${\cal T}$. We say that a symplectic form $\Omega$ on $X$ {\it
tames} ${\cal T}$ if $\Omega(H)>0$ for every subspace $H$ in every ${\cal
T}_{x}$. Then Gromov explains that all the fundamental results about $J$-holomorphic
curves in symplectic manifolds extend to this more general context.

We will now relate this discussion to the rest of this article. Consider
a negative submanifold $T\subset Gr_{2}({\bf R}^{4})$ as before and a map
$f: [0,\infty)\times T \rightarrow \Lambda^{2} {\bf R}^{4}$ with $f(r,\theta)=O(r^{-1})$,
along with all its derivatives. Call the image of the map
$(r,\theta)\mapsto r\theta + f(r,\theta)$ the $f$-deformed cone over $T$. Then we say that a $3$-dimensional submanifold
$P\subset \Lambda^{2}{\bf R}^{4}$ is asymptotic to $T$ if there is a map
$f$ as above such that, outside  compact
subsets, $P$ coincides with the $f$-deformed cone over $T$. This notion immediately
generalises to a pair of constraint manifolds ${\cal P}, {\cal T}$ over a
$4$-manifold.  The prototype
of this picture is to take the Calabi-Yau constraint defined by a volume
form and almost-complex structure, which one readily sees is asymptotic to
the submanifold of complex subspaces. Now suppose that ${\cal P}$ lies in
the positive cone with respect to the wedge-product form, has negative chords
and is asymptotic to ${\cal T}$. Then if $\omega\in {\cal P}_{x}$ and $\theta\in
{\cal T}_{x}$ we have $(\omega- r \theta)^{2} \leq O(r^{-1})$ for large $r$.
Since $\theta^{2}=0$ and $\omega^{2}>0$ we see that $\omega \wedge \theta >0$,
which is the same as saying that $\omega$ tames ${\cal T}$. In other words,
just as in the Calabi-Yau case, a necessary condition for there to be a solution
$\omega \subset {\cal P}$ in a  given cohomology class $h=[\omega]$ is that
there is a taming form in $h$ for ${\cal T}$. It is tempting then to extend
Conjecture 1 to this more general situation.
\begin{conj}
Let $X$ be a compact $4$-manifold and let $\Omega$ be a symplectic form on
$X$. If ${\cal P}$ is a unimodular constraint manifold with negative chords which is
asymptotic to ${\cal T}$, where ${\cal T}$ is tamed by $\Omega$, then there
are $C^{\infty}$ {\it a priori} bounds on a closed form $\omega \subset {\cal
P}$ with $[\omega]=[\Omega]$.
\end{conj}

(There is little hard evidence for the truth of this, so perhaps it is better
considered as a question. Of course, by our results in Section 4 above, it
suffices to obtain an $L^{\infty}$ bound and a modulus of continuity in the
\lq\lq BMO sense'' of Proposition 2.
 Going a very small way in this
direction, it is easy to
show that a taming form for ${\cal T}$ leads to an {\it a priori} $L^{1}$ bound on $\omega$.) 

\subsection{A counterexample}

We can use sophisticated results in symplectic topology to obtain a negative
result--showing that Conjecture 2 cannot be extended to the case where ${\cal
P}$ only has negative tangents.
 
\begin{prop}
There is a simply connected $4$-manifold $X$ with $b^{2}_{+}(X)=1$ and a
$1$-parameter family of  unimodular constraints
${\cal P}_{t} \  (t\in [0,1])$ on $X$  with the following properties
\begin{enumerate}

\item For each $t$, ${\cal P}_{t}$ has negative tangents and  is asymptotic to the manifold ${\cal T}$ associated with an
almost-complex stucture on $X$ which is tamed by a symplectic form $\Omega$.
\item  For $t< 1$ there is a closed $2$-form $\omega_{t}\subset {\cal P}_{t}$ with $[\omega_{t}]\in \langle
[\Omega]\rangle $.
\item The $\omega_{t}$ do not satisfy uniform $C^{\infty}$ bounds as $t\rightarrow
1$. 
\end{enumerate}
\end{prop}
  
In fact we can take the almost-complex structure to be integrable,  $\Omega$
to be a Kahler form , and arrange that, for each parameter value,  ${\cal P}_{t}$  coincides with the standard Calabi-Yau constraint outside a compact
set. 

The proof of Proposition 3 combines some simple general constructions with results of
Seidel. Let $J$ be an almost-complex structure on a $4$-manifold $X$ and
let
$\omega$ be  a symplectic form on $X$. Choose an almost-complex structure $J'$
compatible with $\omega$, and suppose that $J$ and $J'$ are homotopic, through
a 1-parameter family $J(s)$. In other words, for each $s$ we have a 2-dimensional
positive subbundle $\Lambda^{2,0}(s)$. For convenience, suppose that $J(s)=J'$
for $s$ close to $0$ and extend the family to all positive $s$, with $J(s)=J$ for $s\geq 1$. For fixed small $\epsilon$ we define
a subset of $\Lambda^{2}$ by
$$ {\cal P}^{*}= \{ \theta: \theta^{2}=\omega^{2}, \theta \in \Lambda^{2,0}(\epsilon
\vert \theta \vert)\}, $$
where $\vert \theta \vert$ is the norm measured with respect to some arbitrary
metric.  This set ${\cal P}^{*}$ has two connected components, interchanged
by $\theta\mapsto -\theta$, and   it is easy to check that, if $\epsilon $ is sufficiently small, one of these components
 is a constraint manifold ${\cal P}$ with negative tangents, containing $\omega$
and equal to the Calabi-Yau constraint defined by $J$ and $\omega^{2}$ outside
a compact set.

Now we can extend this construction to families. If $\omega_{z}$ is a family
of symplectic forms parametrised by a compact space $Z$, and if the corresponding
homotopy class of maps from $Z$ to the space of almost-complex structures
on $X$ is trivial, then we can construct a family ${\cal P}_{z}$ of unimodular
constraints, equal to the fixed Calabi-Yau constraint outside a compact set
and with ${\cal P}_{z}$ containing $\omega_{z}$. For our application we first
take
$Z$ to be a circle, so we have a loop of symplectic forms $\omega_{z},\ z\in
S^{1}$,  with fixed cohomology class $h=[\omega_{z}]$. Clearly the family ${\cal
P}_{z}$ for $z\in S^{1}$ can be extended over the disc $D^{2}$. So for each
$z\in D^{2}$ we have a ${\cal P}_{z}$ satisfying the conditions of the Proposition.
If the Proposition were false, then it would follow that we could extend the family of symplectic
forms $\omega_{z}$ over the disc, using a simple continuity argument.

The space  ${\cal S}_{\omega}$ of  symplectic forms equivalent to $\omega$ can be
identified with the quotient  of the
identity component of the full diffeomorphism group by the symplectomorphisms
of $(X,\omega)$ which are isotopic to the identity. So the fundamental group
of ${\cal S}_{\omega}$ can be identified with the classes  of symplectomorphisms
isotopic to the identity modulo symplectic isotopy. Now in \cite{kn:S} Seidel gives
examples of Kahler manifolds $(X,\omega)$ satisfying our hypotheses and symplectomorphisms
which are isotopic, but not symplectically isotopic, to the identity. Thus
Seidel's results assert that there are maps from the circle to ${\cal S}_{\omega}$
which cannot be extended to the disc (although the corresponding almost-complex
structures can be). This
conflict with our previous argument completes the proof of Proposition 3.

The author has not yet succeeded in understanding more explicitly the  blow-up behaviour which Proposition 3 asserts must occur.
Seidel's symplectomorphisms are squares of generalised Dehn twists, associated to Lagrangian
$2$-spheres in $(X,\omega)$, and it is tempting to hope that the blow-up
sets should be related to these spheres in some way. 
It should also be noted that Seidel's results depend crucially on fixing
the cohomology class of the symplectic form.



\end{document}